\documentclass{article}
\usepackage{amsmath}
\usepackage{amsfonts}
\usepackage{amssymb}

\begin{document}

\title{Polynomials with small norm on compact Riemannian homogeneous
manifolds}
\author{A. Kushpel \\
Department of Mathematics,\\
University of Leicester}
\date{\today}
\maketitle

\begin{abstract}
Let $H_{k}$, $k \geq 0$ be the sequence of eigenspaces corresponding to the
eigenvalues $0 \leq \theta_{0} \leq \theta_{1} \leq \cdots \leq \theta_{n}
\leq \cdots$ of the Laplace-Beltrami operator $\Delta$ on a compact
Riemannian homogeneous manifold $\mathbb{M}^{d}$ with the normalized
invariant measure $\nu$ and $\mathcal{T}_{n} = \oplus_{k=0}^{n} H_{k}$. We
consider the problem of existence of polynomials $t_{n} \in \mathcal{T}_{n}$
with small norm. Namely, we show that for any $\epsilon \in (0, 1)$ and any
subspace $L_{m} \subset \mathcal{T}_{n}$, $\mathrm{dim} L_{m} \geq \epsilon n
$, there exists such $t_{n} \in L_{m}$ that $\|t_{n}\|_{L_{p}(\mathbb{M}%
^{d}, \nu)} \leq C_{p,q}\|t_{n}\|_{L_{q}(\mathbb{M}^{d}, \nu)}$, where $%
C_{p,q}$ depends just on $p$ and $q$, $1 < q<p < \infty$. In the case $%
p=\infty$ or $q=1$ an extra logarithmic factor appears. This range of
problems has been extensively studied by many authors in the case $\mathbb{M}%
^{d} = \mathbb{T}^{1}$, the unit circle (or compact Abelian group $\mathcal{G%
}$), i.e., when the characters of $\mathcal{G}$ are bounded by $1$. In
general, on compact Riemannian homogeneous manifolds, the eigenfunctions of
the Laplace-Beltrami operator are not uniformly bounded that creates
difficulties of a fundamental nature in applications of known methods and
results. The method, we develop, is based on a geometric inequality between
norms induced by two convex bodies in $\mathbb{R}^{n}$.
\end{abstract}

\textit{MSC:} 41A46, 42B15.\newline

\textit{Keywords:} Homogeneous manifold, volume, L\'evy mean, flat
polynomial.


\section{Introduction}

\label{9619.sec1} The range of problems we consider in this article has been
traditionally studied in the context of random Fourier series and has been
initiated in the classical works of Paley and Zygmund. In many situations it
is difficult or impossible to give explicitly an example of a certain object
having a required property and frequently one gets by with Lebesgue measure
\cite{kahane}, \cite{kahane1}. Problems regarding flat polynomials with
coefficients $\pm 1$ whose uniform norm is close to their $L_{2}(\mathbb{T}%
^{1})$ norm has attracted a lot of attention \cite{borwein}, \cite%
{borwein-erdelyi}. It was shown in \cite{rudin} that for any $N\in \mathbb{N}
$ there is a sequence $\epsilon _{n}=\pm 1$, $1\leq n\leq N$ such that
\[
\left\vert \sum_{k=1}^{N}\epsilon _{n}e^{in\theta }\right\vert <5N^{1/2}.
\]%
This topic has been developed in \cite{kahane2}. It was shown that for all $%
|z|=1$ there is a sequence of polynomials $P(z)=\sum_{m=1}a_{m,n}z^{m}$, $%
|a_{m,n}|=1$ such that
\[
(1-\epsilon _{N})N^{1/2}\leq |P_{N}(z)|\leq (1+\epsilon _{N})N^{1/2},
\]%
where $\epsilon _{N}\leq CN^{-1/17}(\log N)^{1/2}$ as $N\rightarrow \infty $%
. The expected $L_{p}(\mathbb{T}^{1})$ norm of random trigonometric
polynomials $q_{N}(\theta )=\sum_{k=0}^{N}X_{k}e^{ik\theta }$, where $X_{k}$%
, $k\geq 0$ are independent and identically distributed random variables
with mean $0$ and variance $1$ was studied in \cite{borwein1}. In
particular, it was shown that
\[
\frac{\mathbb{E}(\Vert q_{N}\Vert _{p}^{p})}{N^{p/2}}\rightarrow \Gamma
\left( 1+\frac{p}{2}\right) ,\,\,\,N\rightarrow \infty .
\]%
The problem of existence of trigonometric polynomials with special
properties of degree $\leq (M+1)(1+\epsilon )$, $\epsilon >0$, in any
subspace of $L_{2}(\mathbb{T}^{1})$ of codimension $M$ was considered in
\cite{kashin}. It was shown that for any $\epsilon >0$ there is such a
polynomial whose uniform norm is $1$, and such that the sum of the absolute
values of the coefficients is at least $c_{\epsilon }M^{1/2}$.


In this article we consider the problem of existence of polynomials on a
compact homogeneous Riemannian manifold $\mathbb{M}^{d}$ whose $L_{p}(%
\mathbb{M}^{d})$ norm is close to their $L_{2}(\mathbb{M}^{d})$ norm for any
$1 \leq p \leq \infty$ (see Theorem 3).

The method's possibilities are not confined to the theorem proved in the
Section \ref{section4} but can be used in studying more general problems.
The results we derive are apparently new even in the one dimensional case.

\section{Elements of Harmonic Analysis on Compact Riemannian Manifolds}

\label{9619.sec2}

First we give a general definition of function spaces that we consider and
then present various important examples.

{\bf Definition 1.} \label{9619.definition}
\begin{em}
Given a measure space $(\Omega, \nu)$. Let $\Xi =\{\xi_{k}\}_{k
\in \mathbb{N}}$ be a set of orthonormal functions in
$L_{2}(\Omega, \nu)$. Suppose that there exists a sequence
$\{k_{j}\}_{j \in \mathbb{N}}$, $k_1 = 1$, such that for any $j
\in \mathbb{N}$ and some $C>0$
\[
 \sum_{k=k_{j}}^{k_{j+1}-1}
|\xi_{k}(x)|^{2} \leq C d_{j}
\]
a.e. on $\Omega$, where $d_{j} = k_{j+1}-k_{j}$. Then we say that
$(\Omega, \nu, \Xi, \{k_{j}\}_{j \in \mathbb{N}}) \in \cal{K}$.
\end{em}

Let $L_{p}=L_{p}(\Omega, \nu)$ be the usual set of $p$-integrable functions
on $\Omega$. Suppose that $(\Omega, \nu, \Xi, \{k_{j}\}_{j \in \mathbb{N}})
\in \mathcal{K}$. Since all the functions $\xi_{k}$ are a.e. bounded on $%
\Omega$, then for an arbitrary function $\phi \in L_{p}$, $1 \leq p \leq
\infty$ we can construct the sequence $\{c_{k}(\phi)\}_{k \in \mathbb{N}}$,
where $c_{k}(\phi) = \int_{\Omega} \phi \overline{\xi_{k}} d\nu $ and
consider the formal series
\[
\phi \sim \sum_{l=1}^{\infty} \sum_{k_{l}}^{k_{l+1}-1} c_{k}(\phi) \xi_{k}.
\]
The family $\mathcal{K}$ is sufficiently large. We consider compact,
connected, orientable, $d$-dimensional $C^{\infty}$ Riemannian manifold $%
\mathbb{M}^{d}$, with $C^{\infty}$ metric. Let $g$ its metric tensor, $\nu$
its normalized volume element and $\Delta$ its Laplace-Beltrami operator. In
local coordinates $x_l$, $1 \leq l \leq d$,
\[
\Delta = - (\overline{g})^{-1/2} \sum_{k} \frac{\partial}{\partial x_{k}}%
\left(\sum_{j} g^{jk} (\overline{g})^{1/2} \frac{\partial}{\partial x_{j}}%
\right).
\]
Here, $g_{ij} = g(\partial/\partial x_{i}, \partial/\partial x_{j})$, $%
\overline{g} = |\mathrm{det} (g_{ij})|$ and $(g^{ij}) = (g_{ij})^{-1}$. It
is well-known that $\Delta$ is an elliptic, self adjoint, invariant under
isometries, second-order operator. The eigenvalues $\theta_{k}$, $k \geq 0$
of $\Delta$ are discrete, nonnegative and form an increasing sequence $0
\leq \theta_{0} \leq \theta_{1} \leq \cdots \leq \theta_{n} \leq \cdots$
with $+\infty$ the only accumulation point. Corresponding eigenspaces $H_{k}$%
, $k \geq 0$ are finite dimensional, $d_{k} = \dim H_{k} < \infty$, $k \geq 0
$, orthogonal  and  $L_{2}(M^{d}, \nu) = \oplus_{k=0}^{\infty} H_{k}$. Let $%
\{Y_{m}^{k}\}_{m=1}^{d_{k}}$ be an orthonormal basis of $H_{k}$, $H_{k} =
\mathrm{lin} \{Y_{m}^{k}\}_{m=1}^{d_{k}}$.

Recall that a Riemannian manifold $\mathbb{M}^{d}$ is called homogeneous is
its group of isometries ${\mathcal{G}}$ acts transitively on it. Let $H_{j}$%
, $j \geq 0$ be any eigenspace of $\Delta$, $d_{j}= \mathrm{dim} H_{j}$, $%
f_{1}, \cdots f_{d_{j}}$ any orthonormal basis of $H_{j}$, then
\[
\sum_{s=1}^{d_{j}} |f_{s}(x)|^{2} = d_{j},
\]
for any $x \in \mathbb{M}^{d}$ (see, e.g., \cite{9619.gine}). Hence, any
compact, connected, orientable, $d$-dimensional $C^{\infty}$, homogeneous
Riemannian manifold $\mathbb{M}^{d}$, with $C^{\infty}$ metric has the
property ${\mathcal{K}}$. 
Here we give several important examples of such manifolds: \newline
\textbf{1.} A Grassmannian (Grassmann manifold), $\mathbb{G}_{m,n}(\mathbb{R}%
)$ is the space of all m-dimensional subspaces of $\mathbb{R}^{n}$.
Grassmann manifold also appear as coset space $\mathbb{G}_{m,n}(\mathbb{R})
= \mathbf{\mathrm{O}}(n)/\mathbf{\mathrm{O}}(n-m) \times \mathbf{\mathrm{O}}%
(m)$;\newline
\textbf{2.} A complex Grassmannian manifold $\mathbb{G}_{m,n}(\mathbb{C})$
is the space of all $m$-dimensional complex subspaces in $\mathbb{C}^{n}$;
\newline
\textbf{3.} An $n$-torus, $\mathbb{T}^{d}$ is defined as a product of n
circles: $\mathbb{T}^{d} = \mathbb{S}^{1} \times \cdots \mathbb{S}^{1}$. The
$n$-torus can be described as a quotient of $\mathbb{R}^{n}$ under shifts in
any coordinate. That is, the $n$-torus is $\mathbb{R}^{n}$ modulo the action
of the integer lattice $\mathbb{Z}^{n}$ (with the action being taken as
vector addition); \newline
\textbf{4.} The Stiefel manifold, denoted $\mathbb{V}_{k}(\mathbb{R}^{d})$
or $\mathbb{V}_{k,d}$, is the set of all orthonormal $k$-frames in $\mathbb{R%
}^{d}$. That is, it is the set of ordered $k$-tuples of orthonormal vectors
in $\mathbb{R}^{d}$. When $k = 1$, the manifold $\mathbb{V}_{1,d}$ is just
the set of unit vectors in $\mathbb{R}^{d}$; that is, $\mathbb{V}_{1,d}$ is
diffeomorphic to the $d - 1$ sphere, $\mathbb{S}^{d-1}$. At the other
extreme, when $k = d$, the Stiefel manifold $\mathbb{V}_{d,d}$ is the set of
all ordered orthonormal bases for $\mathbb{R}^{d}$. $\mathbb{V}_{d,d}$ is a
principal homogeneous space for $\mathbf{\mathrm{O}}(d)$ and therefore
diffeomorphic to it. In general, the orthogonal group $\mathbf{\mathrm{O}}(d)
$ acts transitively on $\mathbb{V}_{k,d}$ with stabilizer subgroup
isomorphic to $\mathbf{\mathrm{O}}(d-k)$. Therefore $\mathbb{V}_{k,d}$ can
be viewed as the homogeneous space $\mathbb{V}_{k,d} = \mathbf{\mathrm{O}}%
(d)/\mathbf{\mathrm{O}}(d-k). $ \newline
\textbf{5.} The unit complex sphere in $\mathbb{C}^{d}$ is defined as $%
\mathbb{S}_{\mathbb{C}}^{d} = \{z \in \mathbb{C}^{d} |\,\langle z, z \rangle
= 1\}$, where $\langle z, w \rangle = z_{1} \overline{w}_{1} + \cdots z_{d}%
\overline{w}_{d}$, $z,w \in \mathbb{C}^{d}$. \newline
\textbf{6.} A Riemannian manifold is two-point homogeneous if for any set of
four points $x_{1}, y_{1}, x_{2}, y_{2}$ with $d(x_{1}, y_{1}) = d(x_{2},
y_{2})$, $d$ being the Riemannian metric on $\mathbb{M}^{d}$, there exists $%
\phi \in {\mathcal{G}}$ such that $\phi (x_{1}) = x_{2}$ and $\phi(y_{1}) =
y_{2}$. A complete classification of the two-point homogeneous spaces was
given in \cite{9619.wang}.  For information on this classification see,
e.g., \cite{9619.car, 9619.gang, 9619.helgason, 9619.hel1}.  They are: the
spheres $\mathbb{S}^{d}$, $d = 1, 2, 3,\ldots$; the real projective spaces $%
\mathbb{P}^{d}(\mathbb{R})$, $d = 2, 3, 4, \ldots$; the complex projective
spaces $\mathbb{P}^{d}(\mathbb{C})$, $l=d/2$,{\ } $d = 4, 6, 8,\ldots$; the
quaternionic projective spaces $\mathbb{P}^{d}(\mathbb{H})$, $d = 8,
12,\ldots$; the Cayley elliptic plane $P^{16}(\mathrm{Cay})$. The
superscripts here denote the dimension over the reals  of the underlying
manifolds $\mathbb{M}^{d}$.

\section{A geometric inequality}

\label{9619.sec3}


Let $\alpha = (\alpha_1, \cdots, \alpha_n) \in \mathbb{R}^{n}$, $\beta =
(\beta_1, \cdots, \beta_n) \in \mathbb{R}^{n}$ and $\langle \alpha, \beta
\rangle = \sum_{k=1}^{n} \alpha_{k} \beta_{k}$. Let $\|\alpha\|_{(2)} =
\langle \alpha, \alpha \rangle^{1/2}$ be the Euclidean norm on $\mathbb{R}%
^{n}$, $\mathbb{S}^{n-1} = \{\alpha \in \mathbb{R}^{n}:{\ }\|\alpha\|_{(2)}
=1\}$ be the unit sphere in $\mathbb{R}^{n}$, $B_{(2)}^{n} = \{\alpha \in
\mathbb{R}^{n}:{\ }\|\alpha\|_{(2)} \leq 1\}$ be the unit ball in $\mathbb{R}%
^{n}$ and $\mathrm{Vol}_{n}$ be the standard $n$-dimensional volume of
subsets in $\mathbb{R}^{n}$. Let us fix a norm $\|\cdot\|$ on $\mathbb{R}^{n}
$ and denote by $E$ the Banach space $E = (\mathbb{R}^{n}, \|\cdot\|)$ with
the ball $B_{E} = V$. The L\'evy mean $M_{V}$ is defined by
\[
M = M (\mathbb{R}^{n}, \|\cdot\|) = \int_{\mathbb{S}^{n-1}} \|\alpha\|
d\mu(\alpha).
\]
For a convex centrally symmetric body $V \subset \mathbb{R}^{n}$ we define
the polar body $V^{o}$ of $V$ as
\[
V^{o} = \left\{\alpha \in \mathbb{R}^{n}:{\ }\sup_{\beta \in V} |\langle
\alpha, \beta \rangle| \leq 1 \right\}.
\]
The dual space $E^{o} = (\mathbb{R}^{n}, \|\cdot\|_{o})$ is endowed with the
norm
\[
\|\alpha\|_{V^{o}} = \|\alpha\|_{o} = \sup_{\beta \in B_{E}} |\langle
\alpha, \beta \rangle|
\]
and $B_{E^{o}} = V^{o}$.

{\bf Theorem 1.}
\label{Theorem1}
\begin{em}
Let $V$ and $W$ be any convex symmetric bodies in
$\mathbb{R}^{n}$, $V \subset B_{(2)}^{n}$, then for any $n \in
\mathbb{N}$ and $\epsilon > 0$ there is such $0 <\mu_{\epsilon} <
1$ that in any subspace $L_m \subset \mathbb{R}^{n}$, $m=\dim
L_{m} \geq \mu_{\epsilon} n$ there is such $\alpha^{\ast} \in
L_{m}$ that
\[
\frac{\|\alpha^{\ast}\|_{V}}{\|\alpha^{\ast}\|_{W}} \leq
C_{\epsilon} (M_{V})^{1+\epsilon} M_{W^{o}}.
\]
where $C_{\epsilon}$ depends just on $\epsilon$.
\end{em}

{\bf Proof}
From the Urysohn
 inequality (see, e.g., \cite{pisier} p.
6-7) it follows that
\[
\left(\frac{{\rm Vol}_{n}(V)}{{\rm
Vol}_{n}(B_{(2)}^{n})}\right)^{1/n} \leq C_{1}
\int_{\mathbb{S}^{n-1}} \|\alpha\|_{V^{o}} d\mu(\alpha) = C_{1}
M_{V^{o}}
\]
\[
\leq C_{1}
\left(\int_{\mathbb{S}^{n-1}}\|\alpha\|_{V^{o}}^{2}d\mu\right)^{1/2}
= C_{1} M_{V^{o}}
\]
or
\[
{\rm Vol}_{n}(V^{o}) \leq C_{1}^{n} (M_{V})^{n} {\rm Vol}_{n}
(B_{(2)}^{n}).
\]
Comparing the last estimate with the Bourgain-Milman inequality
\cite{bourgain} p. 320,
\[
\left(\frac{{\rm Vol}_{n}(V) \cdot {\rm
Vol}_{n}(V^{o})}{\left({\rm
Vol}_{n}\left(B_{(2)}^{n}\right)\right)^{2}}\right)^{1/n} \geq
C_2,
\]
we get
\begin{equation}  \label{1}
{\rm Vol}_{n} (V)  \geq \left(\frac{C_{2}}{C_{1} M_{V}}\right)^{n}
 {\rm Vol}_{n} (B_{(2)}^{n}).
\end{equation}

Let $W$ be any convex symmetric body in $\mathbb{R}^{n}$. Using
isoperimetric inequality on $\mathbb{S}^{n-1}$ it is possible to
show \cite{bourgain} that for every $0 < \lambda < 1$ there exists
a subspace $L_{m_{1}} \subset \mathbb{R}^{n}$ with $m_{1}={\rm
dim} L_{m_{1}} \geq \lambda n$ such that for any $\alpha \in
L_{m_{1}}$ we have
\begin{equation}  \label{0}
\|\alpha\|_{(2)} \leq C_{3} \frac{M_{W^{o}}}{ (1-\lambda)}
\|\alpha\|_{W}.
\end{equation}
Let $L_{m_{2}} \subset \mathbb{R}^{n}$ be any
$m_{2}$-dimensional subspace. Assume that $m_{1} + m_{2} > n$, so
that $L_{m_{1}} \cap L_{m_{2}} \neq \emptyset$ and
\[
m_3 := \dim (L_{m_{1}} \cap L_{m_{2}}) \geq m_1 + m_2 - n.
\]
Let $(L_{m_{1}} \cap L_{m_{2}})^{\perp}$ be the orthogonal
complement of $L_{m_{1}} \cap L_{m_{2}}$ and $P_{(L_{m_{1}} \cap
L_{m_{2}})^{\perp}}(V)$ be the orthogonal projection of $V$ onto
$(L_{m_{1}} \cap L_{m_{2}})^{\perp}$. Assume that $V \subset
B_{(2)}^{n}$, then
\[
P_{(L_{m_{1}} \cap L_{m_{2}})^{\perp}}(V) \subset P_{(L_{m_{1}}
\cap L_{m_{2}})^{\perp}}\left(B_{(2)}^{m_{3}}\right)
\]
and
\[
{\rm Vol}_{m_{3}} \left(P_{(L_{m_{1}} \cap
L_{m_{2}})^{\perp}}(V)\right) \leq {\rm Vol}_{m_{3}}
\left(P_{(L_{m_{1}} \cap
L_{m_{2}})^{\perp}}\left(B_{(2)}^{m_{3}}\right)\right).
\]
Hence,
\[
{\rm Vol}_{n} (V) = \int_{V} dx = \int_{P_{(L_{m_{1}} \cap
L_{m_{2}})^{\perp}}(V)} {\rm  Vol}_{m_{3}} (V \cap (y + L_{m_{1}}
\cap L_{m_{2}})) dy.
\]
Thus, for any $y \in P_{(L_{m_{1}} \cap L_{m_{2}})^{\perp}}(V)$ by
the Brunn-Minkowski theorem
\[
{\rm Vol}_{m_{3}} (V \cap (y + L_{m_{1}} \cap L_{m_{2}})) \leq
{\rm Vol}_{m_{3}} (V \cap (L_{m_{1}} \cap L_{m_{2}})).
\]
and, therefore,
\[
{\rm Vol}_{n} (V) \leq {\rm Vol}_{m_{3}}(V \cap (L_{m_{1}} \cap
L_{m_{2}})) \cdot {\rm Vol}_{n-m_{3}} \left(P_{(L_{m_{1}} \cap
L_{m_{2}})^{\perp}}(V)\right)
\]
\begin{equation}  \label{2}
\leq {\rm Vol}_{m_{3}}(V \cap (L_{m_{1}} \cap L_{m_{2}})) \cdot
{\rm Vol}_{n-m_{3}} (B_{(2)}^{n-m_{3}})
\end{equation}
Comparing (\ref{1}) and (\ref{2}) we find that for any convex
symmetric body $V \subset B_{(2)}^{n}$ and any $m_{3}$-dimensional
subspace $L_{m_{1}} \cap L_{m_{2}} \subset \mathbb{R}^{n}$,
\begin{equation}  \label{3}
{\rm Vol}_{m_{3}} (V \cap (L_{m_{1}} \cap L_{m_{2}})) \geq
 \left(\frac{C_{2}}{C_{1} M_{V}}\right)^{n}
 \frac{
 {\rm Vol}_{n}(B_{(2)}^{n})}{{\rm Vol}_{n-m_{3}}
 (B_{(2)}^{n-m_{3}}).
 }
\end{equation}
Applying the Santalo inequality (see, e.g. \cite{burago})
\[
\frac{{\rm Vol}_{m_{3}}(V \cap (L_{m_{1}} \cap L_{m_{2}})) \cdot
{\rm Vol}_{m_{3}}((V \cap (L_{m_{1}} \cap L_{m_{2}}))^{o})}{({\rm
Vol}_{m_{3}}(B_{(2)}^{m_{3}}))^{2}} \leq 1
\]
we obtain
\[
{\rm Vol}_{m_{3}}((V \cap (L_{m_{1}} \cap L_{m_{2}}))^{o}) \leq
\frac{({\rm Vol}_{m_{3}}(B_{(2)}^{m_{3}}))^{2}}{{\rm
Vol}_{m_{3}}(V \cap (L_{m_{1}} \cap L_{m_{2}}))}.
\]
Combining this result with the Bieberbach inequality (see, e.g.,
\cite{burago})
\[
2^{m_{3}} {\rm Vol}_{m_{3}} (B_{(2)}^{m_{3}}) ({\rm diam}(V \cap
(L_{m_{1}} \cap L_{m_{2}})))^{-m_{3}}
\]
\[
\leq {\rm Vol}_{m_{3}} ((V \cap (L_{m_{1}} \cap L_{m_{2}}))^{o}),
\]
we get the lower bound for the diameter of the set $V \cap
(L_{m_{1}} \cap L_{m_{2}})$,
\[
{\rm diam} (V \cap (L_{m_{1}} \cap L_{m_{2}})) \geq 2
\left(\frac{{\rm Vol}_{m_{3}}(B_{(2)}^{m_{3}})}{{\rm Vol}_{m_{3}}((V \cap
(L_{m_{1}} \cap L_{m_{2}}))^{o})}\right)^{1/m_{3}}
\]
\begin{equation}  \label{4}
\geq 2 \left(\frac{{\rm Vol}_{m_{3}}(V \cap (L_{m_{1}} \cap
L_{m_{2}}))}{{\rm Vol}_{n} (B_{(2)}^{n})}\right)^{1/m_{3}}.
\end{equation}
Comparing (\ref{3}) and (\ref{4}) we find
\begin{equation}  \label{5}
{\rm diam} (V \cap (L_{m_{1}} \cap L_{m_{2}})) \geq
 \left(\frac{2C_{2}}{C_{1} M_{V}}\right)^{n/m_{3}} \omega_{n, m_{3}},
\end{equation}
where
\[
\omega_{n, m_{3}} := \left(\frac{{\rm Vol}_{n}(B_{(2)}^{n})}{{\rm
Vol}_{n-m_{3}}(B_{(2)}^{n-m_{3}}) \cdot {\rm
Vol}_{m_{3}}(B_{(2)}^{m_{3}}) }\right)^{1/m_{3}}.
\]
Recall that
\[
{\rm Vol}_{n}(B_{(2)}^{n}) = 2 \pi^{n/2} \Gamma \left(\frac{n}{2}+
1\right).
\]
It means that $\omega_{n, m_{3}}$ can be expressed as
\[
\omega_{n, m_{3}} = \left(\frac{\Gamma (n/2+ 1)}{\Gamma((n-m_{3})/
2+ 1) \Gamma(m_{3}/2+ 1)}\right)^{1/m_{3}}.
\]
It is well-known that
\[
\Gamma(z) = z^{z-1/2}e^{-z}(2\pi)^{1/2} \epsilon_{n},{\ }{\
}\lim_{n \rightarrow \infty} \epsilon_{n} = 1,
\]
so that for any $1 \leq m_{3} \leq n$, $n \rightarrow \infty$ we
have
\[
\omega_{n, m_{3}} =
\]
\[
 \left( \frac{
\left(n/2+1\right)^{n/2+1/2}e\epsilon_{n} } { \left((n-m_{3})/2+
1\right)^{(n-m_{3})/2 + 1/2} \left(m_{3}/2+1\right)^{m_{3}/2+1/2}
(2\pi)^{1/2} \epsilon_{n-m_{3}} \epsilon_{m_{3}} }
\right)^{1/m_{3}}
\]
\[
= \left(\frac{(n+2)^{n/2+1/2}}{(n-m_{3}+ 2)^{(n-m_{3})/2 +
1/2}(m_{3}+2)^{m_{3}/2+1/2}}\right)^{1/m_{3}} \left(\frac{2^{1/2}e
\epsilon_{n}}{2\pi^{1/2} \epsilon_{n-m_{3}}
\epsilon_{m_{3}}}\right)^{1/m_{3}}
\]
\[
= \frac{n^{1/2}}{m_{3}^{1/2 + 1/(2m_{3})}}
\]
\[
\times \frac{ \left(1 + 2/n\right)^{n/(2m_{3})+1/(2m_{3})}}
{\left(1 - m_3/n+2/n\right)^{n/(2m_{3})- 1/2 + 1/(2m_{3})}
\left(1+ 2/m_{3}\right)^{1/2 + 1/(2m_{3})} }
\]
\begin{equation}  \label{6}
\times \left(\frac{2^{1/2}e \epsilon_{n}}{2\pi^{1/2}
\epsilon_{n-m_{3}} \epsilon_{m_{3}}}\right)^{1/m_{3}} \asymp
\left(\frac{n}{m_{3}}\right)^{1/2} n^{1/(2m_{3})}.
\end{equation}
Remark that if $m_{3}=\lambda n$, $0 < \lambda < 1$, then
$\omega_{n, m_{3}} \sim (\lambda e)^{-1/2}$ as $n \rightarrow
\infty$.

From (\ref{5}) it follows that for any $L_{m_{2}} \subset
\mathbb{R}^{n}$ there is such $\alpha^{\ast} \in L_{m_{2}}$ that
\begin{equation}  \label{7}
\|\alpha^{\ast}\|_{(2)} \geq \left(\frac{2C_{2}}{C_{1} M_{V}}
\right)^{n/m_{3}} \omega_{n, m_{3}} \|\alpha^{\ast}\|_{V}.
\end{equation}
Recall that $m_{3} = \dim (L_{m_{1}} \cap L_{m_{2}})$. Since
$\alpha^{\ast} \in L_{m_{1}}$ then from (\ref{0}) we get
\begin{equation}  \label{8}
\|\alpha^{\ast}\|_{(2)} \leq C_{3}M_{W^{o}}
\left(\frac{n}{n-m_{1}} \right) \|\alpha^{\ast}\|_{W}.
\end{equation}
Finally, comparing (\ref{7}) and (\ref{8}) we find
\begin{equation}  \label{rex}
\|\alpha^{\ast}\|_{V} \leq
\left(\frac{C_{1}M_{V}}{2C_{2}}\right)^{n/m_{3}}
\frac{C_{3}M_{W^{o}}}{\omega_{n, m_{3}}}\left(\frac{n}{n-m_{1}}
\right)\|\alpha^{\ast}\|_{W}.
\end{equation}
In particular, let $m_{1} = \mu_{1} n$ and $m_{2} = \mu_{2} n$ for
some fixed $\mu_{1} > 0$ and $\mu_{2} > 0$, $1 < \mu_{1} + \mu_{2}
<2$, then from (\ref{6}) and (\ref{rex}) it follows that
\[
\|\alpha^{\ast}\|_{V} \leq C (M_{V})^{1/(\mu_{1}+\mu_{2}-1)}
M_{W^{o}}\|\alpha^{\ast}\|_{W},
\]
where $C > 0$ is an absolute constant.  

\section{Flat Polynomials on $\mathbb{M}^{d}$}

\label{section4}

Let $\Omega$ be a compact space with a normalized measure $\nu$, Fix an
orthonormal system $\Xi = \{\xi_{k}\}_{k \in \mathbb{N}} \subset
L_{2}(\Omega, \nu)$ and a sequence $\{k_{j}\}_{j \in \mathbb{N}}$ such that $%
(\Omega, \nu, \Xi, \{k_{j}\}_{j \in \mathbb{N}}) \in \mathcal{K}$. Let
\[
\Xi^{j} := \mathrm{lin}\{\xi_{k}\}_{k=k_{j}}^{k_{j+1}-1},\,\,\, \Omega_{m}
:= \{j_{1}, \cdots , j_{m}\},\,\,\,\, \Xi(\Omega_{m}) := \mathrm{lin}
\{\Xi^{j_{s}}\}_{s=1}^{m}.
\]
Put $n := \dim \Xi(\Omega_{m}) = \sum_{s=1}^{m} k_{j_{s+1}} - k_{j_{s}} =
\sum_{s=1}^{m} d_{j_{s}}$, where $d_{j_{s}}:= \dim \Xi^{j_{s}}$. Consider
the coordinate isomorphism
\[
J:\,\mathbb{R}^{n} \rightarrow \Xi(\Omega_{m})
\]
that assigns to $\alpha =(\alpha_{1} \cdots , \alpha_{n}) \in \mathbb{R}^{n}$
the function $J\alpha = \xi^{\alpha} = \sum_{l=1}^{n} \alpha_{l} \xi_{j_{l}}
\in \Xi(\Omega_{m})$. Let $X$ be a Banach space such that $\Xi(\Omega_{m})
\subset X$ for any $\Omega_{m} \subset \mathbb{N}$. Put $X_{n} =
\Xi(\Omega_{m}) \cap X$. The definition $\|\alpha\|_{(X_{n})} =
\|\xi^{\alpha}\|_{X}$ induces a norm on $\mathbb{R}^{n}$. Put
\[
B^{n}_{(X_{n})}:=\{\alpha|\,\alpha \in \mathbb{R}^{n},\,\,\|\alpha%
\|_{(X_{n})} \leq 1\},
\]
then $B^{n}_{X_{n}} := J B^{n}_{(X_{n})}$.

A Banach lattice $X$ is $q$-concave, $q < \infty$ if the there is a
constant $C_{q} > 0$ such that
\[
\left(\sum_{i=1}^{n} \|x_{i}(\cdot)\|_{X}^{q}\right)^{1/q} \leq C_{q}(X)
\left\|\left(\sum_{i=1}^{n} |x_{i}(\cdot)|^{q}\right)^{1/q}\right\|_{X}
\]
for any $n \in \mathbb{N}$ and any sequence $\{x_{i}(\cdot)\}_{i=1}^{n}
\subset X$ (see, e.g., \cite{lt}, p. 46).

We will need the following statement.

{\bf Lemma 1.} 
\begin{em}
For any $\xi \in \Xi(\Omega_{m})$, $m \in \mathbb{N}$ we have
\[
\|\xi\|_{L_{p}(\Omega, \nu)}  \leq Cn^{(1/p-1/q)_{+}}
\|\xi\|_{L_{q}(\Omega, \nu)},\,\,\,
\]
where $1 \leq p,q \leq \infty$ and $n := {\rm dim}\,\,
\Xi(\Omega_{m})$.
\end{em}

{\bf Proof}
Consider the function
\[
K_{n}(x,y) := \sum_{\xi_{k} \in \Xi(\Omega_{m})} \xi_{k}(x)
\overline{\xi_{k}(y)}.
\]
Clearly,
\[
K_{n}(x,y) = \int_{\Omega} K_{n}(x,z) K_{n}(z,y) d\nu(z)
\]
and $K_{n}(x,y) = \overline{K_{n}(y,x)}$. Hence,
\[
 \|K_{n}(\cdot,\cdot)\|_{L_{\infty}(\Omega, \nu)} \leq \|K_{n}(y,\cdot)\|_{L_{2}(\Omega, \nu)}
\|K_{n}(x,\cdot)\|_{L_{2}(\Omega, \nu)}
\]
for any $x,y \in \Omega$ and $\|K_{n}(x,\cdot)\|_{L_{2}(\Omega,
\nu)} \leq Cn^{1/2}$, since $(\Omega, \nu, \Xi, \{k_{j}\}_{j \in
\mathbb{N}}) \in \cal{K}$. It means that
\begin{equation}  \label{infty}
 \|K_{n}(\cdot,\cdot)\|_{L_{\infty}(\Omega, \nu)} \leq Cn.\,\,\,\,\,\,
\end{equation}
Let $\xi  \in \Xi(\Omega_{m})$, then applying H\"older inequality
and (\ref{infty}) we get
\[
\|\xi\|_{L_{\infty}(\Omega, \nu)} \leq
\|K_{n}(\cdot,\cdot)\|_{L_{\infty}(\Omega, \nu)}
\|\xi\|_{L_{1}(\Omega, \nu)},
\]
or
\[
\|I\|_{L_{1}(\Omega, \nu) \cap
\Xi(\Omega_{m}) \rightarrow L_{\infty}(\Omega, \nu) \cap
\Xi(\Omega_{m})} \leq Cn,
\]
 where
$I$ is an embedding operator. Trivially, $\|I\|_{L_{p}(\Omega,
\nu)} \cap \Xi(\Omega_{m}) \rightarrow L_{p}(\Omega, \nu) \cap
\Xi(\Omega_{m}) =1$, $1 \leq p \leq \infty$. Hence, using
Riesz-Thorin interpolation Theorem and an embedding arguments for
any $\xi \in \Xi(\Omega_{m})$ we obtain
\[
\|\xi\|_{L_{p}(\Omega,
\nu)}  \leq Cn^{(1/p-1/q)_{+}} \|\xi\|_{L_{q}(\Omega,
\nu)},\,\,\,1 \leq p,q \leq \infty,
\]
where
\[
(a)_{+} := \left\{
\begin{array}{cc}
  a, & a \geq 0, \\
  0, & a <0. \\
\end{array}
\right.
\]

In the case $\Omega_{m} = \{1, \cdots ,m\}$ the estimates of respective
L\'evy means have been obtained in \cite{ku3}. Using Lemma 1 we
can generalize our result to an arbitrary index set $\Omega_{m}$.

{\bf Theorem 2.} 
\begin{em}
Let $(\Omega, \nu, \Xi, \{k_{l}\}) \in {\mathcal
K}$ and $X$ is a $2$-concave, then for an arbitrary $\Omega_{m}$,
\begin{equation}  \label{star}
M(\mathbb{R}^{n}, \|\cdot\|_{(X_{n})}) \leq C_{X},\,\,\,X_{n} = X \cap
\Xi(\Omega_{m})
\end{equation}
where $n := {\rm dim} \Xi(\Omega_{m})$ and $C_{X} >0$ is
independent on $n \in \mathbb{N}$. In particular,
\begin{equation}  \label{star1}
M(\mathbb{R}^{n}, \|\cdot\|_{(L_{p}(\Omega, \nu) \cap \Xi(\Omega_{m}))})
\leq C \left\{
\begin{array}{cc}
  p^{1/2}, & p < \infty, \\
 (\log n)^{1/2}, & p = \infty, \\
\end{array}
\right.
\end{equation}
where  $C>0$ is an absolute constant.
\end{em}

Remark that different estimates of L\'evy means have been obtained in \cite%
{ku-lev1} - \cite{9619.ku2}. We are prepared now to prove main result of
this article.

{\bf Theorem 3.} 
\begin{em}
Assume that $\max\{M_{j^{-1} \circ (X \cap \Xi_{n})}, M_{j^{-1}
\circ (Y \cap \Xi_{n})^{o}}\} < C$ for any $n \in \mathbb{N}$ and
some absolute constant $C>0$ and $B_{X} \subset
B_{L_{2}(\mathbb{M}^{d})}$. Then in any subspace $J^{-1} \circ
L_{s} \subset \Xi(\Omega_{m})$ there exists such polynomial
$t_{n}^{\ast}$ that
\begin{equation}  \label{24}
\|t_{n}^{\ast}\|_{X} \leq C_{X,Y}\|t_{n}^{\ast}\|_{Y}.
\end{equation}
In particular, the inequality (\ref{24}) is valid if  $(\Omega,
\nu, \Xi, \{k_{j}\}_{j \in \mathbb{N}}) \in {\cal K}$, $X$ is
$2$-concave and  $\|\cdot\|_{(B_{Y} \cap \Xi_{n})^{o}} \leq
\|\cdot\|_{J^{-1} \circ (Y_{1} \cap \Xi_{n})}$ for some 2-concave
$Y_{1}$ and any $n \in \mathbb{N}$. Let $(\Omega, \nu, \Xi,
\{k_{j}\}_{j \in \mathbb{N}}) \in {\cal K}$, $X=L_{p}(\Omega,
\nu)$, $Y=L_{q}(\Omega, \nu)$, $1 \leq p,q \leq \infty$, then for
an arbitrary spectrum $\Omega_{m}$, $n = \dim \Xi(\Omega_{m})$ and
any subspace $J^{-1} \circ L_{s} \subset \Xi(\Omega_{m})$ there
exists a polynomial $t_{n}^{\ast} \in J^{-1} \circ L_{s}$ such
that
\[
\|t_{n}^{\ast}\|_{L_{p}(\Omega, \nu)} \leq C \varrho_{n}
\|t_{n}^{\ast}\|_{L_{q}(\Omega, \nu)},
\]
where
\[
\varrho_{n} = \left\{
\begin{array}{cc}
  1, & 1 < q,p < \infty, \\
  (\log n)^{1/2}, & 1 \leq q \leq p < \infty, \\
(\log n)^{1/2}, & 1 < q \leq p \leq \infty, \\
\log n, & 1 \leq q \leq p \leq \infty \\
\end{array}
\right.
\]
and $C>0$ is an absolute constant.
\end{em}

{\bf Proof}
Applying Theorem 2 and Theorem 1 for a fixed $\epsilon \in (0, 1)$
and the inequality
\[
\|\cdot\|_{(B_{Y} \cap \Xi_{n})^{o}} \leq
\|\cdot\|_{Y_{1} \cap \Xi_{n}},
\]
where $Y_{1}$ is a 2-concave, we
get
\[
\frac{\|\alpha^{\ast}\|_{V}}{\|\alpha^{\ast}\|_{W}} \leq C
M_{V}^{1+\epsilon} M_{W^{o}} \leq C_{X, Y_{1}}
\]
Hence, using (\ref{star}), for any $L_{s} \subset \mathbb{R}^{n}$,
$s = \dim L_{s} \geq \epsilon n$, $\epsilon \in (0, 1)$, one can
find such $\alpha^{\ast} \in L_{s}$ that
\[
\|\alpha^{\ast}\|_{V} \leq C_{X, Y_{1}} \|\alpha^{\ast}\|_{W}.
\]
It means that in any subspace $J \circ L_{s}$ there exists such
$t_{n}^{\ast} \in J \circ L_{s} \subset \Xi(\Omega_{m})$ that
\[
\|t_{n}^{\ast}\|_{X} \leq C_{X, Y_{1}} \|t_{n}^{\ast}\|_{Y}.
\]
In the case  $X=L_{p}(\Omega, \nu)$ and $X=L_{q}(\Omega, \nu)$ we
use (\ref{star1}) to get a similar estimate
\[
\|t_{n}^{\ast}\|_{L_{p}(\Omega, \nu)} \leq C \varrho_{n}
\|t_{n}^{\ast}\|_{L_{q}(\Omega, \nu)},
\]
where $C>0$ is an absolute constant. 



\begin{thebibliography}{99}
\bibitem{bourgain} Bourgain, J. K., Milman, V. D. (1987). New volume ratio
properties for convex bodies in $\mathbb{R}^{n}$. \textit{Invent. Math.},
\textbf{88}, p. 319--340.

\bibitem{borwein-erdelyi} Borwein, P., Erd\'elyi, (1996). Questions about
polynomials with $\{0, -1, +1\}$ coefficients, \textit{Constructive
Approximation}, \textbf{12}, (3), 439-442.

\bibitem{borwein1} Borwein, P., Lockhart, R., (2001). The expected $L_{p}$
norm of random polynomials, \textit{Proc. Amer. Math. Soc.}, \textbf{129},
(5), 1463-1472.

\bibitem{borwein} Borvein, P., (2002). \textit{Paul Erd\"os and polynomials}%
, Bolyai Soc. Math. Stud., 161-174.

\bibitem{burago} Burago, Yu. D., Zalgaller, V. A. (1980). \textit{Geometric
Inequalities}, Springler-Verlag, Berlin, Aidelberg, New York, London, Paris,
Tokyo.


\bibitem{9619.car} Cartan, E. (1929). Sur la determination d'un systeme
orthogonal complet dans un espace de Riemann symetrique clos, \textit{%
Circolo Matematico di Palermo, Rendiconti}, \textbf{53}, 217--252.





\bibitem{9619.gang} Gangolli, R. (1967). Positive definite kernels on
homogeneous spaces and certain stochastic processes related to L\'{e}vy's
Browian motion of several parameters, \textit{Ann. Inst. H. Poincar\'{e}},
\textbf{3}, 121--226.

\bibitem{9619.gine} Gin\'{e}, E. M. (1975). The Addition Formula for the
Eigenfunctions of the Laplacian, \textit{Advances in Mathematics}, \textbf{18%
}, 102--107.


\bibitem{9619.helgason} Helgason, S. (1962). \emph{Differential Geometry and
Symmetric Spaces}, Academic Press, New York.

\bibitem{9619.hel1} Helgason, S. (1965). The Radon transform on Euclidean
spaces, compact two-point homogeneous spaces and Grassmann manifolds,
\textit{\ Acta Math.}, \textbf{113}, 153--180.



\bibitem{kahane} Kahane, J. -P., (1968). \textit{Some random series of
functions}, Nauka, Moscow.

\bibitem{kahane1} Kahane, J. -P. (1971). The technique of using random
measures and random sets in harmonic analysis, \textit{Advances in
probability and related topics}, Dekker, New York, \textbf{1}, 65-101.

\bibitem{kahane2} Kahane, J. -P. (1980). Sur les polyn\^omes \`acoefficients
unimodulaires, \textit{Bull. London Math. Soc.}, Dekker, New York, \textbf{12%
}, (5), 321.

\bibitem{kashin} Kashin, B. S. (1980) Some properties of the space of
trigonometric polynomials with a uniform norm, \textit{Trudy Mat. Inst.
Steklov.}, \textbf{145}, (250), 111--116.

\bibitem{ku-lev1} Kushpel, A. K. (1992).  Estimates of L\'evy means and
medians of some  distributions on a sphere, In: \emph{Fourier Series and
Their Applications }, Kiev, Inst. Math., 49--53.

\bibitem{ku-lev2} Kushpel, A. K. (1993) Estimates of the Bernstein Widths
and Their Analogs, \emph{Ukrainian Mathematical Journal,\newline
Plenum Publishing, New York }, \textbf{45}, 1, 54--69.

\bibitem{9619.ku1} Kushpel, A. K. (1999). L\'evy Means associated with
Two-Point Homogeneous Spaces and Applications, \textit{49 Semin\'ario
Brasileiro de An\'alise}, Campinas, SP, 807--823.

\bibitem{9619.ku2} Kushpel, A. K. (1999). Estimates of $n$-Widths and $%
\epsilon$-Entropy of Sobolev's Sets on Compact Globally Symmetric Spaces of
Rank 1, \textit{50 Semin\'ario Brasileiro de An\'alise}, S\~ao Paulo, SP,
53--66.

\bibitem{ku3} Kushpel, A. K., Tozoni, S. A. (2007) On the Problem of Optimal
Reconstruction, \textit{Journal of Fourier Analysis and Applications},
\textbf{13}, (4), 459--475.


\bibitem{lt} Lidenstrauss, J., Tzafriri, L (1979). \textit{Classical Banach
Spaces II, Function Spaces}, Springer-Verlag.






\bibitem{pisier} Pisier, G. (1989). \textit{The volume of Convex Bodies and
Banach Space Geometry,} Cambridge: Cambridge University Press.


\bibitem{rudin} Rudin, W., (1980). Sur les polyn\^omes \`a coefficients
unimodulares, \textit{Bull. London Math. Soc.}, \textbf{12}, (5) 321-342.


\bibitem{9619.wang} Wang, H. C. (1952). Two-point homogeneous spaces,
\textit{Ann. of Math.}, \textbf{55}, 177--191.



\end{thebibliography}
\end{document}